\documentclass[11pt]{article}
\usepackage{amsmath}
\usepackage{amssymb}
\usepackage{gastex}
\newtheorem{lemma}{{\sc Lemma}}[section]
\newtheorem{theo}{{\sc Theorem}}

\newtheorem{defn}{{\em Definition}}
\newtheorem{rem}{{\em Remark}}

\newtheorem{algo}{{\sc Algorithm}}

\begin{document}

\title{Hamiltonian Cycles in Polyhedral Maps}
\author{{Dipendu Maity  and Ashish Kumar Upadhyay}\\[2mm]
{\normalsize Department of Mathematics}\\{\normalsize Indian Institute of Technology Patna }\\{\normalsize Patliputra Colony, Patna 800\,013,  India.}\\
{\small \{dipendumaity, upadhyay\}@iitp.ac.in}}
%\date{May 2009}
\maketitle

\vspace{-5mm}

\hrule

\begin{abstract} We present a necessary and sufficient condition for existence of a contractible, non-separating and noncontractible separating Hamiltonian cycle in the edge graph of polyhedral maps on surfaces. In particular, we show the existence of contractible Hamiltonian cycle in equivelar triangulated maps. We also present an algorithm to construct such cycles whenever it exists.
\end{abstract}

{\small

{\bf AMS classification\,:} 57M20, 57N05, 05C38.

{\bf Keywords\,:} Contractible Hamiltonian Cycles, Non-separating Hamiltonian Cycles,  Noncontractible Separating Hamiltonian Cycles, Proper Graphs in Polyhedral Maps.
}

\bigskip

\hrule

\section{Introduction and Definitions}

By a graph $G := (V, E)$ we mean a finite simple graph with vertex set $V$ and edge set $E$. Let $G_{1}(V_{1}, E_{1})$ and $G_{2}(V_{2}, E_{2})$ be two simple graphs embedded on a surface. Then $G_{1}\cup G_{2}$ is a graph $G(V, E)$ where $V = V_{1}\cup V_{2}$ and $E = E_{1}\cup E_{2}$. Similarly, $G_{1}\cap G_{2}$ is a graph $G(V, E)$ where $V = V_{1}\cap V_{2}$ and $E = E_{1}\cap E_{2}$, see \cite{BondyMurthy}. A {\em connected component} or {\em component} of a graph is a  subgraph in which any two vertices are connected by a path. We call $H_1$, $H_2$, \ldots, $H_k$ components of $G$ if $G = H_1 \cup H_2 \cup \dots \cup H_k$ and there is no path between the vertices of $H_i$ and $H_j$ for $i \neq j$. A {\em surface} $S$ is a $2$-dimensional manifold which is connected, compact and without boundary. A {\em map} on a surface $S$ is an embedding of a graph $G$ such that the closure of each component of $S \setminus G$ is a $p$-gonal $2$-disc for $p \geq 3$. The components are called {\em faces} of the map and the vertices and edges of the embedded graph $G$ are called vertices and edges of the map. The map $M$ is called a {\em polyhedral map} if intersection of any two faces of $M$ is either empty, a vertex or an edge, see \cite{brehm_schulte}. A map $M$ is said to be a {\em triangulation} of the surface if each face of the map is a $3$-gon. We call $G$ the edge graph of $M$ and denote it by $EG(M)$. A map is called {\em equivelar triangulation} if each vertex of the map $M$ has same degree. We will use the terms map and polyhedral map interchangeably to mean a polyhedral map. A path $P$ in a graph $G$ is a subgraph of $G$, such that the vertex set of $P$ is $V(P) = \{v_1, v_2, \dots, v_n\}\subseteq V(G)$ and $v_{i}v_{i + 1}$ are edges in $P$ for $1 \leq i \leq n - 1$. A path $P$ in $G$ is said to be a cycle if $v_n v_1$ is also an edge in $P$ and denote it by $C(v_1, v_2, \dots, v_n)$. A graph without any cycles is called an {\em acyclic} graph, see \cite{BondyMurthy}. Length $l(C)$ of a cycle $C$ is the number of edges in $C$. See \cite{mohar thomassen} for details about graphs on surfaces and \cite{BondyMurthy} for graph theory related terminology.

A cycle of a graph is said to be a Hamiltonian cycle if it contains all the vertices of the graph. In this article, we are interested in finding out whether a Hamiltonian cycle exists in the edge graph of a polyhedral map? In this context, Tutte \cite{Tutte} has shown that every $4$-connected planar graph has a Hamiltonian cycle. In 1970, Gr$\ddot{u}$nbaum \cite{grunbaum} has conjectured that every $4$-connected graph which admits an embedding in the torus has a Hamiltonian cycle. We summarize the known partial results related to the solution of Gr$\ddot{u}$nbaum's conjecture. Altshuler \cite{a1} has shown that every $4$ and $6$-connected equivelar map of types $\{4, 4\}$ and $\{3, 6\}$, respectively, on the torus has Hamiltonian cycles. Brunet and Richter\cite{brunet richter} have shown that every $5$-connected triangulations on torus is Hamiltonian and then, Thomas and Yu\cite{thomas yu} improved this result for any $5$-connected graph on the torus. Brunet, Nakamoto and Negami \cite{brunet_nakamoto_negami} have shown that every $5$-connected triangulated Klein bottle is Hamiltonian. In \cite{fuji_naka_ozeki} it is shown that a $3$-connected bipartite graph embeddable in torus has a Hamiltonian cycle if it is balanced and each vertex of one of its partite sets has degree four. Kawarabayashi and Ozeki\cite{kawarabayashi ozeki} have shown that every $4$-connected triangulated torus is Hamiltonian. In \cite{maity upadhyay1} we have extended this result to any $4$-connected semi-equivelar map, see definition of {\em semi-equivelar map} in \cite{maity upadhyay1}. In this article, we show every equivelar triangulated map is Hamiltonian. Topologically, such cycles in maps may or may not be homotopic to the generators of fundamental group of the surface on which they lie. The cycles which are homotopic to a generator are called {\em essential or non-separating cycles}. Those which are homotopic to a point are called {\em inessential or contractible cycles}. Those which are not homotopic to both point and generator are called {\em noncontractible separating cycles}. Among these, contractible Hamiltonian cycles have been investigated in \cite{maity upadhyay0, Upadhyay}. We have given a necessary and sufficient condition for existence of contractible Hamiltonian cycles by introducing {\em proper tree} in the dual map (defined later) of polyhedral maps. Separating cycles in triangulations of the double torus studied in \cite{jennings}. Archdeacon\cite {Archdeacon} has given a survey of such cycles in maps. Thus naturally led to think about existence of a graph in the dual which may give information about such types of Hamiltonian cycles. In this article our focus is on finding such cycles in polyhedral maps. We accomplish this by introducing {\em proper graph} (see, Definition \ref{defn1}).

We begin with some definitions which will be needed in the course of proof of main Theorem \ref{thm1}, \ref{thm2}, \ref{thm3} and \ref{thm4}. For more details on these topics one may also refer to \cite{mcmullen schulte}.

Let $v$ be a vertex of a map $K$. The {\em degree} of a vertex $v$ is denoted by $\deg(v)$ and it is the number of edges incident with $v$. If we denote the number of vertices, edges and faces of $K$ by $f_0(K)$, $f_1(K)$ and $f_2(K)$ respectively, then the {\em Euler characteristic of $K$} is the integer $\chi(K) = f_0(K) - f_1(K) + f_2(K)$. The {\em dual map} of a map $M$ is by definition a map on same surface $S$ as $M$ which has a vertex corresponding to each face of $M$, and an edge joining two neighboring faces for each edge in $M$. Let $K$ be a polyhedral map on a surface $S$. Let $M$ denote the dual map of $K$ and $C(u_1,$ \dots, $u_r)$ be a non-separating cycle in $K$. We consider the dual edges of $C$ and put in $E$. Let $V = \cup_{e \in E}V(e)$. The graph $G := (V, E)$ is said to be {\em dual graph} corresponding to the cycle $C$ in $M$. In Section \ref{ex} we give examples of maps on torus. These are well known maps of type $\{3, 6\}$ and $\{6, 3\}$ and are examples of mutually dual maps. Let $K'$ be a subset of set of faces of $K$ and $D = \displaystyle\cup_{\small{\sigma \in K'}}\sigma$. If $D$ is topologically a 2-disc then we will call it a {\em 2-disc} in $K$. A vertex $v$ of $EG(M)$ is called {\em cut vertex} if its removal graph becomes disconnected. Similarly, we call a cycle $C$ in $M$ is {\em cut cycle} if it divides $F_M$ into disjoint set of faces. Similarly, {\em cut graph} is defined. Consider a polyhedral map $K$ on a surface $S$ that has $n$ vertices. Let $M$ denote the dual map of $K$. Let $T := (V, E)$ denote a tree in the edge graph $EG(M)$ of $M$. We say that $T$ is a {\em proper tree}\cite{maity upadhyay0} if the following conditions hold : $(1)$ $\displaystyle\sum_{i = 1}^{ k }\deg(v_{i}) = n + 2(k-1)$, where $V = \{v_{1},v_{2},\dots,v_{k}\}$ and $\deg(v)$ denotes degree of $v$ in $EG(M)$ $(2)$ whenever two vertices $u_{1}$ and $u_{2}$ of $T$ lie on a face $F$ in $M$, a path $P[u_{1},u_{2}]$ joining $u_{1}$ and $u_{2}$ in the boundary $\partial{F}$ of $F$ is a subtree of $T$, and $(3)$ any path $P$ in $T$ which lies in a face $F$ of $M$ is of length at most $q - 2$, where $q =$ length of $(\partial{F})$.

Let $M$ be a polyhedral map which is the dual map of a polyhedral map $K$ on $n$ vertices. We call $(V_M, E_M, F_M)$ the map where $V_M$, $E_M$ and $F_M$ denote the set of vertices, edges and faces of $M$.

\begin{defn}\label{defn1} Let $G = (V, E)$ denote a subgraph of $EG(M) : = (V_M, E_M)$. We say that $G$ is an {\em admissible graph} if the following conditions hold\,:
\begin{enumerate}
\item $|E| = n$, that is, $G$ consists of $n$ edges.
\item For each $F \in F_M$, the set $F \cap E$ contains exactly two edges of $E$.
\item There exists a finite sequence of $n$ faces, namely, $F_1$, $F_2$, \ldots, $F_n$ such that $F_i \cap F_{i+1} \in E$ for $1 \leq i \leq n-1$ and $F_1 \cap F_n \in E$.
\end{enumerate}

We say that $G$ is a  proper graph of type-I if the graph $H = (V_M, E_M\setminus E)$ is connected and call $G$ a proper graph of type-II if $H$ has exactly two components one of which is a proper tree. If the graph $G(V_M, E_M\setminus E)$ has two components where none of them is proper tree then $G$ is called a proper graph of type-III.

\end{defn}

The main result of this article is\,:

\begin{theo}\label{thm1} The edge graph $EG(M)$ of a  map $M$ has a non-separating Hamiltonian cycle if and only if the edge graph of the dual map of $M$ has a proper graph of type-I.
\end{theo}

In \cite{maity upadhyay0} we have presented a necessary and sufficient condition for existence of a contractible Hamiltonian cycle in a polyhedral map using a tree. Here we present the criterion in terms of a more general graph and show\,:

\begin{theo}\label{thm2} The edge graph $EG(M)$ of a  map $M$ has a contractible Hamiltonian cycle if and only if the edge graph of the dual map of $M$ has a proper graph of type-II.
\end{theo}

\begin{theo}\label{thm3} The edge graph $EG(K)$ of a  map $K$ has a noncontractible separating Hamiltonian cycle if and only if the edge graph of corresponding dual map of $K$ has a proper graph of type-III.
\end{theo}

We have defined the graphs, namely, proper graph of type-I, proper graph of type-II and proper graph of type-III. A triangulated map is called {\em equivelar} if degree of each vertex is same. We show the existence of proper graph of type-II in equivelar triangulated maps. This leads to the result,\

\begin{theo}\label{thm4} Equivelar triangulated map is Hamiltonian.
\end{theo}

In Section \ref{ex}, we give examples of proper graphs of type-I, II, III and their corresponding cycles. In Section \ref{type1}, \ref{type2} and \ref{type3} we present some properties of proper graphs and proceed to prove the main result of this article. In Section \ref{chctm}, we show the existence of proper graph of type-II in equivelar triangulated maps. We give two different algorithms in Section \ref{algorithm} and \ref{algorithm1} to detect a non-contractible Hamiltonian cycle. The algorithm \ref{algorithm} is not implemented by us where as the algorithm \ref{algorithm1} is implemented.

\section {Example of proper graphs and their corresponding cycles}\label{ex}

We consider two well known maps of types $\{3, 6\}$ and $\{6, 3\}$ on the torus. Let $M_1 := \{u_{11}u_{12}u_{14},$ $u_{11}u_{13}u_{14},$ $ u_{12}u_{13}u_{15},$ $u_{12}u_{14}u_{15},$ $u_{13}u_{14}u_{16},$ $u_{13}u_{15}u_{16},$ $u_{14}u_{15}u_{17},$ $u_{14}u_{16}u_{17},$ $ u_{15}u_{16}u_{11},$ $u_{15}u_{17}u_{11},$ $u_{16}u_{17}u_{12},$ $u_{16}u_{11}u_{12},$ $u_{17}u_{11}u_{13},$ $u_{17}u_{12}u_{13}\}$ and $K_1 := \{[v_{1},$ $ v_{2},$ $ v_{3},$ $ v_{8},$ $ v_{7},$ $ v_{6}],$ $[v_{3},$ $ v_{4},$ $ v_{5},$ $ v_{10},$ $ v_{9},$ $ v_{8}],$ $[v_{5},$ $ v_{6},$ $ v_{7},$ $ v_{12},$ $ v_{11}, $ $v_{10}],$ $ [v_{7},$ $ v_{8},$ $ v_{9},$ $ v_{14},$ $ v_{13},$ $ v_{12}],$ $[v_{9},$ $ v_{10},$ $ v_{11},$ $ v_{2},$ $ v_{1},$ $ v_{14}],$ $[v_{11},$ $ v_{12},$ $ v_{13},$ $ v_{4},$ $ v_{3}, $ $v_{2}],$ $ [v_{13},$ $ v_{14},$ $ v_{1},$ $ v_{6},$ $ v_{5},$ $ v_{4}] \}$. The $M_1$  and $K_1$ are dual to each other. We consider proper graphs of type-I, II in $K_1$. Let $G_1=(V_1,$ $E_1)$ where $V_1$ = $\{v_{i}$ $|$ $i \in \{1,$ $\dots,$ $ 14\}\}$ and $E_1$ = $\{v_{i}v_{5+i}$ $|$ $i \in \{1,$ $3,$ $5,$ $7,$ $9\}\}$ $\cup$ $\{v_{i}v_{9+i}$ $|$ $i \in \{2,$ $4\}\}$ and  $G_2 = (V_2,$ $E_2)$ where $V_2$ = $\{v_{i}$ $|$ $i \in \{1,$ $2,$ $3,$ $4,$ $5,$ $8,$ $10,$ $11,$ $12,$ $13,$ $14\}\}$ and $E_2$ = $\{ v_{1}v_{2},$ $v_{3}v_{8},$ $v_{4}v_{5},$ $v_{10}v_{11},$ $v_{11}v_{12},$ $v_{12}v_{13},$ $v_{13}v_{14} \}$. The graph $G_1$ is of type-I and $G_2$ is of type-II. Now consider dual of $G_1$ and $G_2$ in $M_1$. We get the cycle $C_1(u_{11},$ $u_{12},$ $\dots,$ $u_{17})$ which is corresponding to $G_1$ and non-contractible. Also, the cycle $C_2(u_{11},$ $u_{14},$ $u_{15},$ $u_{13},$ $u_{17},$ $u_{12},$ $u_{16})$ which is corresponding to disconnected graph $G_2$ and contractible.

Let $K : = \{126,$ $267,$ $237,$ $378,$ $348,$ $489,$ $459,$ $59a,$ $56a,$ $6ab,$ $67b,$ $7bc,$ $78c,$ $18c,$ $189,$ $129,$ $29a,$ $23a,$ $3ab,$ $34b,$ $4bc,$ $45c,$ $15c,$ $16l,$ $6lf,$ $ldf,$ $dfg,$ $deg,$ $egh,$ $5eh,$ $5hi,$ $56i,$ $6ij,$ $6fj,$ $fjk,$ $fgk,$ $1gk,$ $1gh,$ $1lh,$ $lhi,$ $ldi,$ $dij,$ $dej,$ $ejk,$ $5ek,$ $15k\}$ be a triangulation of double torus. Let $M :=\{[a_1,$ $b_7,$ $b_6,$ $b_5,$ $c_5,$ $f_1,$ $e_1,$ $e_2,$ $e_3,$ $c_6],$ $[a_1,$ $a_2,$ $a_3,$ $b_9,$ $b_8,$ $b_7],$ $[a_3,$ $a_4,$ $a_5,$ $c_2,$ $c_1,$ $b_9],$ $[a_5,$ $a_6,$ $a_7,$ $c_4,$ $c_3,$ $c_2],$ $ [a_7,$ $a_8,$ $a_9,$ $d_5,$ $d_4,$ $d_3,$ $e_9,$ $f_1,$ $c_5,$ $c_4],$ $[a_9,$ $b_1,$ $b_2,$ $a_2,$ $a_1,$ $c_6,$ $c_7,$ $d_7,$ $d_6,$ $d_5],$ $[b_2,$ $b_3,$ $b_4,$ $a_4,$ $a_3,$ $a_2],$ $[b_4,$ $b_5,$ $b_6,$ $a_6,$ $a_5,$ $a_4],$ $[b_6,$ $b_7,$ $b_8,$ $a_8,$ $a_7,$ $a_6],$ $[b_8,$ $b_9,$ $ c_1,$ $b_1,$ $a_9,$ $a_8],$ $[c_1,$ $c_2,$ $c_3,$ $b_3,$ $b_2,$ $b_1],$ $[c_3,$ $c_4,$ $c_5,$ $b_5,$ $b_4,$ $b_3],$ $[c_6,$ $c_7,$ $c_8,$ $e_5,$ $e_4,$ $e_3],$ $[c_8,$ $c_9,$ $d_1,$ $e_7,$ $e_6,$ $e_5],$ $[d_1,$ $d_2,$ $d_3,$ $e_9,$ $e_8,$ $e_7],$ $[d_7,$ $d_8,$ $d_9,$ $c_9,$ $c_8,$ $c_7],$ $[d_9,$ $e_1,$ $e_2,$ $d_2,$ $d_1,$ $c_9],$ $[e_2,$ $e_3,$ $e_4,$ $d_4,$ $d_3,$ $d_2],$ $[e_4,$ $e_5,$ $e_6,$ $d_6,$ $d_5,$ $d_4],$ $[e_6,$ $e_7,$ $e_8,$ $d_8,$ $d_7,$ $d_6],$ $[e_8,$ $e_9,$ $f_1,$ $e_1,$ $d_9,$ $d_8]\}$. The map $M$ is dual of $K$. Let $E_3 = \{a_1c_6,$ $a_9b_1,$ $b_1c_1,$ $ b_2b_3,$ $b_3b_4,$ $c_4c_5,$ $e_9f_1,$ $e_8d_8,$ $d_7d_8,$ $c_7c_8,$ $e_4e_5,$ $e_5e_6,$ $d_1e_7,$ $d_2d_3,$ $d_2e_2,$ $e_1e_2\}$ and $V_3 := \{a_1,$ $c_6,$ $a_9,$ $b_1,$ $c_1,$ $b_2,$ $b_3,$ $b_4,$ $c_4,$ $c_5,$ $e_9,$ $f_1,$ $e_8,$ $d_8,$ $d_7,$ $c_7,$ $c_8,$ $e_4,$ $e_5,$  $e_6,$ $d_1,$ $e_7,$ $d_2,$ $d_3,$ $d_2,$ $e_2,$ $e_1\}$. The graph $G_3(V_3, E_3)$ satisfies all the properties of proper graph of type-III. So, graph $G_3$ is of type-III. Consider dual of $G_3$ in $K$. We get the cycle $C_3(1,$ $6,$ $a,$ $b,$ $7,$ $c,$ $5,$ $k,$ $j,$ $f,$ $l,$ $i,$ $d,$ $e,$ $h,$ $g)$. The cycle $C_3$ is noncontractible separating Hamiltonian cycle.

\section{Properties of proper graph of type-I\,:}\label{type1}

Let $K$ be a polyhedral map on a surface $S$ with $n$ vertices and $M$ denote its dual. Let $C(u_1,$ \dots, $u_n)$ be a non-separating Hamiltonian cycle in $K$. We consider $G(V, E)$ the dual graph corresponding to the cycle $C$ in $M$. We claim the following lemma.

\begin{lemma}\label{lem1} The graph $G$ contains $n$ edges.
\end{lemma}

\noindent{\sc Proof of Lemma}\ref{lem1} : The cycle $C$ consists of $n$ edges. By the definition of duality, we get an edge in $M$ for an edge of $C$. So, we get $n$ distinct edges corresponding to $n$ edges of $C$. Therefore, $G$ contains $n$ edges. This completes the proof lemma \ref{lem1}.
%The graph $G$ is proper graph in a polyhedral map $M$. Suppose $G$ contains a cycle, say $C(u_1, \dots, u_r)$. Consider dual edges of the edges of cycle $C$. We get $r ( \ge 3)$ edges corresponding to $r$ edges of $C$. If $C$ bounds a face then there will be at least three edges which have common vertex.
\hfill$\Box$

\begin{lemma}\label{lem2} Let $F$ be a face of $M$. Then, $F \cap E$ contains exactly two edges of $EG(M)$.
\end{lemma}

\noindent{\sc Proof of Lemma}\ref{lem2} : Let $v$ be a vertex of $K$. The link of $v$ is a cycle, say $C_1(v_1,$ $\dots,$ $v_r)$. Since the cycle $C$ is Hamiltonian cycle in $EG(K)$, so, it contains all the vertices of $K$ and hence, $v \in V(C)$. Also, the cycle $C$ passes through two vertices, say $v_i,$ $v_j$ of $C_1$ where the edges $vv_i$ and $vv_j$ are in $C$. Since the degree of each vertex in $C$ is 2 therefore, there is no edge $vv_k$ for which $v_k$ $\in$ $V(C_1)\setminus \{v_i,$ $v_j\}$ in $C$. Now consider the dual face corresponding to $v$ in $M$ and denote it by $F$. Also, consider the edges which are dual of $vv_i$ and $vv_j$. We get that exactly two edges are common in $G$ and $F$ and that there is no other edge of $G$ which also belongs to $F$. Since, choice of $v$ was arbitrary this is true for all the faces of $M$. Therefore, any face of $M$ contains exactly two edges of $G$. This completes the proof of lemma \ref{lem2}.
\hfill$\Box$

\begin{lemma}\label{lem3} There exists a face sequence $F_1,$ $F_2,$ $\dots,$ $F_n$ in $M$ such that $F_i$ $\cap$ $F_{i+1}$ $\in$ $E$ for $1$ $\leq$ $i$ $\leq$ $n-1$ and $F_1$ $\cap$ $F_n$ $\in$ $E$.
\end{lemma}

\noindent{\sc Proof of Lemma}\ref{lem3} : Let $F_i$ denote the faces in $M$ which are the dual faces corresponding to  vertices  $u_i$ of $C(u_1, \dots, u_n)$ in $K$. Then $F_i$ and $F_{i+1}$ have an edge in common for $1 \leq i \leq n-1$ and so does $F_1$ and $F_n$. Since $u_{i-1}$ and $u_{i+1}$ lie in the link of $u_i$, by the argument in proof of previous lemma \ref{lem2} it follows that the common edges lie in $E$.  Hence we get a sequence of faces $F_1,$ $\dots,$ $F_n$ where $F_1$ $\cap$ $F_n$ $\in E$ and  $F_i$ $\cap$ $F_{i+1}$ $\in E$ for $1$ $\le$ $i$ $\le n-1$. This completes the proof of lemma \ref{lem3}.
\hfill$\Box$

\begin{lemma}\label{lem4} The graph $G(V_M,$ $E_M\setminus E)$ is connected.
\end{lemma}

\noindent{\sc Proof of Lemma}\ref{lem4} : Let $u_{t_1}$ and $u_{t_r}$ be two vertices of $G_1 := G(V_M, E_M\setminus E)$. We show the existence of a path between $u_{t_1}$ and $u_{t_r}$ in $G_1$. We call $FS(A, B) : = \{A = F_1, F_2, \dots, F_r = B ~ | ~ F_i \cap F_{i+1} \text{ is an edge for } 1 \le i \le r-1\}$ a face sequence between two faces $A$ and $B$. Let $CEFS(A, B) := \{ e_i ~ | ~ e_i = F_i \cap F_{i+1} \text{ and } F_i, F_{i+1} \in FS(A,B)\}$. The set $CEFS(A, B)$ is a set of common edges between the successive faces of $FS(A, B)$. Let $F_{u_{t_1}}$ and $F_{u_{t_r}}$ be two dual faces corresponding to $u_{t_1}$ and $u_{t_r}$ respectively in $K$. The map $K$ is connected polyhedral map. So, we get a face sequence $FS(F_{u_{t_1}}, F_{u_{t_r}})$. We claim that there exists a face sequence $FS(F_{u_{t_1}}, F_{u_{t_r}})$ such that $CEFS(F_{u_{t_1}}, F_{u_{t_r}}) \cap E(C) = \emptyset$. Suppose, there is no such face sequence $FS(F_{u_{t_1}}, F_{u_{t_r}})$ in $K$. That is, for every $FS(F_{u_{t_1}}, F_{u_{t_r}})$, $CEFS(F_{u_{t_1}}, F_{u_{t_r}}) \cap E(C) \not= \emptyset$. This imply that, the cycle $C$ is cut cycle. Let $C$ divides $F(K)$ into disjoint set of faces, namely, $H_1$, \dots, $H_t$. Then, each edge of $|H_i|$ is either identified or on the boundary. The edges which are on the boundary of $|H_i|$ belongs to $E(C)$. So, $\cup_{i=1}^{t}\partial |H_i| = C$. Consider adjacent faces of $C$. It has two sequences, namely, $F_1, F_2, \dots, F_{r1}$ and $Q_1, Q_2, \dots, Q_{r2}$ such that non empty intersection $F_i \cap Q_j$ is part of $C$. Without loss of generality, we assume that $F_i \in H_1$ and $Q_j \in H_2$ for all $i$ and $j$. Then, $\partial |H_i| = C$ for $ i =  1, 2$. This gives, $t = 2$. Thus, $H_1 \cup H_2 =  F(K)$ and $\partial |H_1| \cap \partial |H_2| = C$. So, the cycle $C$ is separating. This is a contradiction as $C$ is non-separating. Therefore, there exits a face sequence $FS(F_{u_{t_1}}, F_{u_{t_r}})$ with $CEFS(F_{u_{t_1}}, F_{u_{t_r}}) \cap E = \emptyset$. We consider dual graph, say, $G_d$ of the face sequence $FS(F_{u_{t_1}}, F_{u_{t_r}})$. The graph $G_d$ is connected as between any two consecutive faces of $FS(F_{u_{t_1}}, F_{u_{t_r}})$ have a common edge. So, the graph $G_d$ contains a path$(u_{t_1}\rightarrow u_{t_r})$ in $G(V_M, E_M\setminus E)$. This is true for any two vertices of $G(V_M, E_M\setminus E)$. Therefore, the graph $G(V_M, E_M\setminus E)$ is connected. This completes the proof of lemma \ref{lem4}.
\hfill$\Box$

\begin{lemma}\label{lem5} Let $G_1$ be a proper graph of type-I in $M$ and $C_1$ be its dual in $K$. Suppose, $G(V_M, E_M\setminus E(G_1))$ is connected. Then, $C_1$ is non-separating.
\end{lemma}

\noindent{\sc Proof of Lemma}\ref{lem5} : The graph $G(V_M, E_M\setminus E(G_1))$ is connected. We show the cycle $C_1$ is non-separating. Suppose $C_1$ is separating. We use some arguments and definitions of lemma \ref{lem4}. The cycle $C_1$ divides $F(K)$ into two sets, say $N_1$ and $N_2$ where $N_1 \cap N_2 = \emptyset$, $N_1 \cup N_2 = F(K)$ and $|N_1| \cap |N_2| = C_1$. Let $F_i \in N_1$ and $F_j \in N_2$ be two faces of $K$ and $v_k$ be dual vertices corresponding to $F_k$. Consider a face sequence $FS(F_i, F_j)$ between $F_i$ and $F_j$ in $K$. By assumption, $F_i \in N_1$, $F_j \in N_2$ and $|N_1| \cap |N_2| = C_1$. So, $CEFS(F_i, F_j) \cap E(C_1) \not=\emptyset$ and it does for any face sequence $FS(F_i, F_j)$. Let $CEFS(F_i, F_j) \cap E(C_1) = \{e_1, e_2, \dots, e_s\}$. Consider dual of the face sequence $FS(F_i, F_j)$ and denote it by $H$. The graph $H$ is connected. The vertices $v_i, v_j \in H$. So, there is a path $P(v_i\rightarrow v_j)$ in $H$. Let $d_i$ be dual edges corresponding to $e_i$. So, path $P(v_i\rightarrow v_j)$ contains some $d_i$. Thus, for arbitrary path $P(v_i \rightarrow v_j)$ in $H$, $E(P) \cap E(G_1) \not= \emptyset$. We restrict the path $P$ in $G(V_M, E_M\setminus E(G_1))$. We get a disconnection between $v_i$ and $v_j$. So, there is no path between $v_i$ and $v_j$ in $G(V_M, E_M\setminus E(G_1))$. This gives a contradiction as $G(V_M, E_M\setminus E(G_1))$ is connected. Therefore $C_1$ is non-separating. This completes the proof of lemma \ref{lem5}.
\hfill$\Box$

%\begin{lemma}\label{lem5} Let $C_1$ be a cycle in $K$ and $G_1$ denote its dual in $M$. If $G(V_M, E_M\setminus E(G_1))$ is a connected graph then $C_1$ is non-separating.
%\end{lemma}
%
%\noindent{\sc Proof of Lemma}\ref{lem5} : The graph $G(V_M, E_M\setminus E(G_1))$ is connected. We show that $C_1$ is non-separating. Suppose $C_1$ is separating. Then the cycle $C_1$ bounds a $2$-disk, say $D_1$. Let $F_i$ and $F_o$ be two faces of $K$ where $F_i \in F(D_1)$ and $F_o \not\in F(D_1)$ where $F(D_1)$ denotes the faces of $D_1$. We consider a face sequence between $F_i$ and $F_o$ in $K$ where any two consecutive faces have a common edge. Since $F_i \in F(D_1)$ and $F_o \not\in D_1$, so, for arbitrary face sequence between $F_i$ and $F_o$ there exists some pair of consecutive faces which have common edges in $E(C_1)$. Let $F_i, \dots, F_o$ be a face sequence between $F_i$ and $F_o$. We denote the common edges in $E(C_1)$ of the face sequence by $e_1, e_2, \dots, e_t$. We consider dual of the face sequence. This is a connected graph which contains a path $P(v_o \rightarrow v_i)$. The path $P$ contains some of the dual edges of $e_1, \dots, e_t$. This gives, $E(P) \cap E(G_1) \not= \emptyset$. Therefore, if we restrict the path $P$ in $G(V_M, E_M\setminus E(G_1))$ we get a disconnection between $v_o$ and $v_i$. So, there is no path between $v_o$ and $v_i$. This gives a contradiction as $G(V_M, E_M\setminus E(G_1))$ is connected. Therefore $C_1$ is non-contractible. This completes the proof of lemma \ref{lem5}.
%\hfill$\Box$

\section{Properties of proper graph of type-II\,:}\label{type2}

Let $K$ be a polyhedral map on $n$ vertices. We consider the dual map of $K$ and denote it by $M$. Let $C(u_1,$ \dots, $u_n)$ be a contractible Hamiltonian cycle in $K$ and $G := (V, E)$ be its dual graph. We claim the following lemma.

\begin{lemma}\label{lem6} The graph $G(V_M,$ $E_M\setminus$ $E)$ has two components where one of them is proper tree.
\end{lemma}

\noindent{\sc Proof of Lemma}\ref{lem6} : The cycle $C$ is contractible. So, cycle $C$ bounds a $2$-disk. Let $D \subset F(K)$ such that $\partial |D| = C$ and $|D|$ is homeomorphic to $2$-disk. Let $K_1 = D$ and $K_2 = F(K) \setminus D$. Then, $\partial |K_1| = C$, $\partial |K_2| = C$ and $\partial |K_1| \cap \partial|K_2| = C$. Let $M_i$ be dual of $K_{i}$. Then, $EG(M) = EG(M_1)\cup EG(M_2)\cup G$. Consider adjacent faces of cycle $C(u_1, u_2, \dots, u_n)$. It has two sequences, namely, $F_1, F_2, \dots, F_n$ and $Q_1, Q_2, \dots, Q_n$ where $F_i \cap Q_i = u_iu_{i+1}$ for all $ i \in \{1, \dots, n-1\}$ and $F_n \cap Q_n = u_nu_{1}$. Without loss of generality, we assume that $F_i \in K_1$ and $Q_i \in K_2$ for all $i$. For each $i$, consider dual edges corresponding to $u_iu_{i+1}$ and denote it by $v_iv_{i+1}$. By definition of $G$, $E = \{v_iv_{i+1} ~ | ~ 1 \le i \le n-1\} \cup \{v_nv_1\}$. Let $F_u \in K_1$ and $F_w \in K_2$ and $u,$ $w$ be dual vertices corresponding to faces $F_u$ and $F_w$ respectively. By construction, $\partial |K_1| \cap \partial|K_2| = C$. So, $Q(u \rightarrow w) \cap E \not= \emptyset$. It is true for arbitrary path $Q(u \rightarrow w)$. Thus, the set $E$ is cut set. The set $E$ divides $EG(M)$ into two components $EG(M_1)$ and $EG(M_2)$ and $G(V_M, E_M \setminus E) = EG(M_1) \cup EG(M_2)$. So, $G(V_M, E_M \setminus E)$ has two components $EG(M_{1})$ and $EG(M_{2})$. Again, by assumption, $D$ is $2$-disk, $\partial |D| = C$ and $M_1$ is dual of $D$. So, $EG(M_{1})$ is proper tree (by the $lemma~5.2$\cite{maity upadhyay0}). Therefore, one of the component of $G(V_M,$ $E_M\setminus$ $E)$ is proper tree. This completes the proof of the lemma \ref{lem6}.
\hfill$\Box$

%Since $C$ is boundary cycle of $|K_i|$ for $i \in \{1, 2\}$, so, in dual of $K_i$, we do not consider the dual edges $E$ of $C$. So, the dual maps of $K_1$ and $K_2$ are two connected components and denoted by $M_{1}$ and $M_{2}$ respectively. Hence, any path between the vertices of $M_{1}$ and $M_{2}$ must contains the edges of $E$. Therefore, there is no path$(u \rightarrow v)$ in $G(V_M, E_M \setminus E)$ for $u \in V_{M_{1}}$ and  $v \in V_{M_{2}}$. The map $M_1$ is $2$-disk as $|K_1| = |D|$. We denote the edge graph of $M_1$ by $G_{M_{1}}$. Similarly, $|K_2|$ is a map with a boundary, so, $M_2$ is a map with a boundary. We denote its edge graph by $G_{M_{2}}$. So, the graph $G(V_M,$ $E_M\setminus E) = G_{M_{1}}$ $\cup$ $G_{M_{2}}$. It has two components $G_{M_{1}}$ and $G_{M_{2}}$. Again, since $D$ is $2$-disk and bounded by Hamiltonian cycle $C$ and $M_1$ is dual of $D$, so, the edge graph $G_{M_{1}}$ of $M_1$ is a proper tree. Hence, one of the component of $G(V_M,$ $E_M\setminus$ $E)$ is a proper tree. This completes the proof of the lemma \ref{lem6}.

\begin{lemma}\label{lem7} Let $H_1$ be a proper graph of type-II in $M$ and $C_1$ be its dual in $K$. Suppose, graph $G(V_M, E_M\setminus E(H_1))$ consists two components, say $G_1$ and $G_2$. The graph $G_1$ is a proper tree. Then, $C_1$ is contractible.
\end{lemma}

\noindent{\sc Proof of Lemma}\ref{lem7} : The graph $G_1$ is proper tree. By {\em lemma 5.2}\cite{maity upadhyay0}, the dual of $G_1$ is $2$-disk, say $D$ and bounded by a Hamiltonian cycle, say $C_2$. By remark \ref{rem0}, $G_1$ defines $H_1$ and $H_1$ is proper graph of type-II. By the argument of theorem \ref{thm5}, $H_1$ and $G_1$ define same cycle. That is, $\partial |D| = C_1 =C_2$. Therefore, the cycle $C_1$ is contractible Hamiltonian cycle. This completes the proof of lemma \ref{lem7}.
\hfill$\Box$

\section{Properties of proper graph of type-III\,:}\label{type3}

Let $K$ be a polyhedral map on $n$ vertices. We consider the dual map of $K$ and denote it by $M$. Let $C(u_1, \dots, u_n)$ be a noncontractible separating Hamiltonian cycle in $K$ and $G := (V, E)$ be its dual graph in $M$. We claim the following lemma.

\begin{lemma}\label{lem8} The graph $G(V_M,$ $E_M\setminus$ $E)$ has two components where none of them is proper tree.
\end{lemma}

\noindent{\sc Proof of Lemma}\ref{lem8} : The cycle $C$ is separating. So, $C$ divides $F(K)$ into two sets, say $T_1$ and $T_2$ where $T_1 \cap T_2 = \emptyset$, $T_1 \cup T_2 = F(K)$ and $\partial |T_1| \cap \partial |T_2| = C$. We follow similar argument of lemma \ref{lem6}. Consider adjacent faces of $C$. It has two sequences, namely, $F_1, F_2, \dots, F_n$ in $T_1$ and $Q_1, Q_2, \dots, Q_n$ in $T_2$ where $F_i \cap Q_i = u_iu_{i+1}$ for $ 1 \le i \le n-1$ and $F_n \cap Q_n = u_nu_{1}$. For each $i$, consider dual edges of $u_iu_{i+1}$ and denote it by $v_iv_{i+1}$. By definition of $G$, $E = \{v_iv_{i+1} : 1 \le i \le n-1\} \cup \{v_1v_n\}$. Consider dual maps corresponding to $T_i$ and denote it by $X_i$. Thus, $EG(M) = EG(X_1) \cup EG(X_2) \cup G$ and $G$ is a cut graph. Therefore, $G(V_M, E_M \setminus E) = EG(X_1) \cup EG(X_2)$ has two components $EG(X_{1})$ and $EG(X_{2})$. Now, we show that none of $EG(M_{1})$ and $EG(M_{2})$ is proper tree. Suppose $EG(M_{1})$ is proper tree. Then, by remark \ref{rem0}, the graph $G$ is proper graph of type-II. This gives the cycle is contractible. This is a contradiction as the cycle is non-contractible. Therefore, the  graph $G(V_M,$ $E_M\setminus$ $E)$ has two components where none of them is proper tree. This completes the proof of the lemma \ref{lem8}.
\hfill$\Box$

\begin{lemma}\label{lem9} Let $H_1$ be a proper graph of type-III in $M$ and $C_1$ be its dual in $K$. Suppose, graph $G(V_M, E_M\setminus E(H_1))$ consists two components, say $G_1$ and $G_2$. None of $G_1$ and $G_2$ is proper tree. Then $C_1$ is a noncontractible separating cycle in $K$.
\end{lemma}

\noindent{\sc Proof of Lemma}\ref{lem9} : The graph $G(V_M, E_M\setminus E(H_1))$ has two components $G_1$ and $G_2$. That is, $G(V_M, E_M\setminus E(H_1)) = G_1 \cup G_2$ and $G_1 \cap G_2 = \emptyset$. We show that $C_1$ is non-contractible. Suppose $C_1$ is contractible. Then, the cycle $C_1$ bounds $2$-disk $|D|$. Consider dual of $D$ which is one component of $G(V_M, E_M\setminus E(H_1))$. Let dual of $D$ be $G_1$. By remark \ref{rem0}, the graph $G_1$ is a proper tree. This is a contradiction as none of $G_1$ and $G_2$ is proper tree. So, the cycle $C_1$ is non-contractible. The graph $G(V_M, E_M\setminus E(H_1))$ is disconnected. So, by lemma \ref{lem4}, the cycle $C_1$ separating. Therefore, the cycle $C_1$ is noncontractible separating cycle. This completes the proof of the lemma \ref{lem9}.
\hfill$\Box$

\section{Proof of theorem \ref{thm1}\,:}

%\noindent {\sc Proof of Theorem } \ref{thm1} :

Let $K$ be a polyhedral map on $n$ vertices and $M$ denote its dual map. Suppose, $C(u_1, \dots, u_n)$ is non-separating Hamiltonian cycle. Let $G := (V, E)$ be dual graph corresponding to $C$ in $M$. By lemma \ref{lem1}, the graph $G$ consists of $n$ edges. That is, $|E| = n$. By lemma \ref{lem2}, $\# E(G) \cap F =2~ \forall ~ F \in F_M$ ($\#S$ denotes the number of elements). By lemma \ref{lem3}, there exists a face sequence, namely, $F_1, F_2, \dots, F_n$ where $F_1 \cap F_n \in E$ and  $F_i \cap F_{i+1} \in E$ for $1 \le i \le n-1$. By lemma \ref{lem4}, $G(V_M, E_M\setminus E)$ is connected. Therefore, the graph $G$ is a proper graph of type-I.

Let $G = (V, E)$ be a proper graph of type-I. By the third property of $G$, there exists a face sequence, say $F_1, F_2, \dots, F_n$ of $M$ where $F_1 \cap F_n \in E$ and  $F_i \cap F_{i+1} \in E$ for $1 \le i \le n-1$. Let $u_i$ be dual vertex corresponding to face $F_i$ in $K$. By definition, $F_i \cap F_{i+1}$ is an edge in $G$ which is common edge between the faces $F_i$ and $F_{i+1}$. So, the edge $F_i \cap F_{i+1}$ is dual of $u_iu_{i+1}$. Put, $u_iu_{i+1}$ in a set, say $E_1$. The sets $F_i \cap F_{i+1}$ for $i \in \{1, 2, \dots, n-1\}$ contain exactly one edge. So, the set $E_1$ contains exactly $n$ edges. Thus, $E_1 = \{u_1u_2,$ $u_2u_3,$\dots, $u_iu_{i+1},$\dots, $u_{n-1}u_n,$ $u_nu_1\}$. Define $C :=C(u_1,$ $u_2,$ \dots, $u_i,$ $u_{i+1},$ \dots, $u_n)$. The cycle $C$ contains all the vertices of $K$. So, $C$ is Hamiltonian. Again, $G(V_M, E_M\setminus E)$ is connected. By lemma \ref{lem5}, $C$ is non-separating. Therefore, $EG(K)$ contains non-separating Hamiltonian cycles. This completes the proof of theorem \ref{thm1}.
\hfill$\Box$

\section{Proof of theorem \ref{thm2}\,:}

%\noindent {\sc Proof of Theorem } \ref{thm2} :

Let $K$ be a polyhedral map on $n$ vertices and $M$ denote its dual. Suppose $C(u_1, \dots, u_n)$ is a contractible Hamiltonian cycle. Let $G := (V, E)$ be dual graph corresponding to $C$ in $M$. We use similar argument of theorem \ref{thm1}. By lemma \ref{lem1}, \ref{lem2} and \ref{lem3}, $G$ is an admissible graph. By lemma \ref{lem6}, the graph $G(V_M, E_M\setminus E)$ has two components and one of them is proper tree. Therefore, the graph $G$ is a proper graph of type-II.

Let $G = (V, E)$ be a proper graph of type-II. We follow similar argument of theorem \ref{thm1}. By the third property of $G$, there exists a sequence of faces, say $F_1, F_2, \dots, F_n$ in $M$ where $F_1 \cap F_n \in E$ and  $F_i \cap F_{i+1} \in E$ for $1 \le i \le n-1$. Let $u_i$ be dual vertex corresponding to face $F_i$. This gives an edge set $\{u_1u_2,$ $u_2u_3,$\dots, $u_iu_{i+1},$\dots, $u_{n-1}u_n,$ $u_nu_1\}$ in $K$. Define $C :=C(u_1,$ $u_2,$ $\dots,$ $u_i,$ $u_{i+1},$ $\dots,$ $u_n)$. The cycle $C$ contains all the vertices of $K$. So, $C$ is Hamiltonian. Graph $G(V_M, E_M\setminus E)$ is of type-II. So, it has two components and one of them is proper tree. By lemma \ref{lem7}, $C$ is contractible. Therefore, $EG(K)$ contains contractible Hamiltonian cycles. This completes the proof of theorem \ref{thm2}.
\hfill$\Box$

\section{Proof of theorem \ref{thm3}\,:}

Let $K$ be a polyhedral map on $n$ vertices and $M$ denote its dual. Suppose $C(u_1, \dots, u_n)$ is a noncontractible separating Hamiltonian cycle. Let $G := (V, E)$ be dual graph corresponding to $C$ in $M$. We use similar argument of theorem \ref{thm1}. By lemma \ref{lem1}, \ref{lem2} and \ref{lem3}, graph $G$ is admissible. By lemma \ref{lem8}, the graph $G(V_M, E_M\setminus E)$ has two components and none of them is proper tree. Therefore, the graph $G$ is a proper graph of type-III.

Let $G = (V, E)$ be a proper graph of type-III. We follow similar argument of theorem \ref{thm1} to show the dual of $G$ is a Hamiltonian cycle. By the third property of $G$, there exists a sequence of faces, say $F_1, F_2, \dots, F_n$ of $M$ where $F_1 \cap F_n \in E$ and  $F_i \cap F_{i+1} \in E$ for $1 \le i \le n-1$. Let $u_i$ be dual vertex corresponding to face $F_i$. This gives an edge set $\{u_1u_2,$ $u_2u_3,$\dots, $u_iu_{i+1},$\dots, $u_{n-1}u_n,$ $u_nu_1\}$ in $K$. Define $C :=C(u_1,$ $u_2,$ $\dots,$ $u_i,$ $u_{i+1},$ $\dots,$ $u_n)$. The cycle $C$ contains all the vertices of $K$. So, $C$ is Hamiltonian. Graph $G(V_M, E_M\setminus E)$ has two components and none of them is proper tree. So by lemma \ref{lem9}, $C$ is noncontractible and separating. Therefore, $EG(K)$ contains noncontractible separating Hamiltonian cycles. This completes the proof of theorem \ref{thm3}.
\hfill$\Box$

\section{Properties of proper tree in equivelar triangulation \,:}

Let $K$ be a triangulated map on $n$ vertices and $M$ be its dual. Then we claim that,\

\begin{lemma}\label{lem8}  Let $T$ be a tree in $M$. Tree $T$ satisfies the following three properties of proper tree.
\begin{enumerate}
\item Whenever two vertices $u_{1}$ and $u_{2}$ of $T$ lie on a face $F$ in $M$, a path $P[u_{1},u_{2}]$ joining $u_{1}$ and $u_{2}$ in the boundary $\partial{F}$ of $F$ is a subtree of $T$,
\item any path $P$ in $T$ which lies in a face $F$ of $M$ is of length at most $q - 2$, where $q =$ length of $(\partial{F})$, and
\item $T$ touches all the faces of $M$.
\end{enumerate}

Then, $\#V(T)=n-2$.
\end{lemma}

\noindent{\sc Proof of Lemma}\ref{lem8} : Let $\#V(T) = m $ $<$ $(n-2)$. Consider dual of $T$ in $K$. We get a $2$-disk $D_T$ which consists of $m$ triangles. So, the cycle $\partial D_T$ contains $m+2$ number of vertices. By assumption, $m < (n-2)$. That is,  $(m+2) <n$. So, we get a vertex $u$ in $K$ which is not belongs to $\partial D_T$. Take dual of $u$ in $M$ and denote it by $F_u$. Since, $u \not\in \partial D_T$, so $T$ does not touch $F_u$. This a contradiction as $T$ touches all the faces of $M$. Hence, $\#V(T) = m $ $\ge$ $(n-2)$. Suppose, $\#V(T) = m $ $>$ $(n-2)$. Similarly, consider dual of $T$ in $K$. Since, $T$ satisfies $2^{nd}$ and $3^{rd}$ property of proper tree, so we get a $2$-disk $D_T$, see \cite{Upadhyay}. The $2$-disk $D_T$ consists of $m$ triangles. So, the cycle $\partial D_T$ contains $m+2$ number of vertices. That is, $m+2 > n$ as $m > n-2$. But, map $K$ consists of $n$ vertices. So, there is a repetition of vertices on the boundary of $D_T$. Thus, $D_T$ is not a $2$-disk. This is a contradiction. Therefore, $m = n-2$.
\hfill$\Box$

\section{Proof of theorem \ref{thm4} \,:}\label{chctm}

Let $K$ be a $q$-equivelar triangulated map on $n$ vertices and $M$ be its dual. The map $M$ is of type $\{q, 3\}$. Cut $M$ along the generators. It gives a planer representation of $M$ and denote it by $N_1$. Planer representation $N_1$ consists of only $q$-gons. We redefine the planer representation $N_1$ as follows : let $F_0$ be a face in $N_1$. There are exactly $q$ number of adjacent faces and each face has an edge which is intersection with the face $F_0$. Consider these $q$ faces and identify with $F_0$. We denote it by $L$. Delete faces of $L$ from $N_1$. Put unidentified edges of $L$ in a set, say $H_1$ which are not repeated. Then, consider faces $F_i$ in $N_1$ such that $F_i \cap H_1 \not=\emptyset$ and identify with the faces of $L$. Delete those $F_i$ from $N_1$. Continue until edges in $L$ either identified or repeated exactly twice. Thus, we get a planer representation from $N_1$ and denote it by $N$. Now, we construct a proper tree $T = (V, E)$ in $M$. Consider all the repeated edges of $N$ and put in a set, say $B$. Let $F(u_1, \dots, u_q)$ be a face of $N$ where $E(F) \cap B = \emptyset$. The face $F(u_1, \dots, u_q)$ denotes a $q$-gon which is bounded by the cycle $C(u_1, \dots, u_q)$. The face $F$ exists as $M$ is a polyhedral map. Consider path $P(u_1\rightarrow \dots \rightarrow u_{q-1}$). Length of $P$ is $q-2$. Put, $V = \{u_1, u_2, \dots, u_{q-1}\}$, $E = \{u_1u_2, u_2u_3, \dots, u_{q-2}u_{q-1}\}$ and $D = \{F\}\cup \{F_l \in F(N)| F\cap F_l \not=\emptyset\}$. The vertex $u_q \not\in V$ and edges $u_{q-1}u_q$, $u_qu_1\not\in E$. There are two faces at $u_q$ except $F$, say $F_1$ and $F_2$. Let $F_i(u_{q-1}, u_q, w_2, w_3, \dots, w_{q-1}) \cap B =\emptyset$. Here, $\# F_i \cap V = 1$. The cycle $\partial F_i$ consists of two paths $P_{i1}(w_{q-2}\rightarrow w_{q-1}\rightarrow u_{q-1} \rightarrow u_q \rightarrow w_2 \rightarrow w_3)$ and $P_{i2}(w_3 \rightarrow w_4 \rightarrow \dots \rightarrow w_{q-2})$ such that $P_{i1} \cup P_{i2} = \partial F_i$ and all adjacent faces of $P_{i1}\setminus\{w_{q-2}, w_3\}$ are in $D$. Consider adjacent faces of $P_{i2}$ in $F(N)\setminus D$ and put in $D$. We also put the edges of $P_{i2}$ and edges $u_{q-1}w_{q-1}, w_{q-1}w_{q-2}$ in $E$ and vertices $w_3, w_4, \dots, w_{q-1}$ in $V$. We continue with this process for the faces $F_x \in D$ such that $\# F_x \cap V = 1 $ and $F_x \cap B = \emptyset$. Let $F_j$ be a face such that $F_j \cap E \not= \emptyset$ and $F_j \cap B = \emptyset$. The cycle $\partial F_j$ consists of two paths $P_{j1}$ and $P_{j2}$ such that $\partial F_j = P_{j1} \cup P_{j2}$, $E(P_{j1}) \cap E(P_{j2}) = \emptyset$ and $P_{j1}$ is part of tree $T$. Let $z, w$ be two end vertices of $P_{j2}$ and $P_{j1}$, $P_{j1}(w \rightarrow x_1 \rightarrow \dots \rightarrow x_r\rightarrow z)$ and $P_{j2}(w \rightarrow y_1 \rightarrow \dots \rightarrow y_s\rightarrow z)$. Consider adjacent faces of $P_{j2}$ in $F(N)\setminus D$ and put in $D$. We put the vertices of $P_{j2}\setminus \{z, y_s\}$ in $V$ and the edges of $P_{j2}\setminus \{z, y_s\}$ in $E$. We continue with the above process for the faces $F_x \in D$ such that $ F_x \cap E \not= \emptyset $ and $F_x \cap B = \emptyset$. There are three faces at each vertex and any two faces have a common edge. All the faces have same length. So, at each step, we choose a face $F_0$ in $D$. There are the following two possibilities. Either all adjacent faces of $F_0$ are in $D$ or some are not in $D$. Suppose, some adjacent faces of $F_0$ are not in $D$. We have the following three possibilities. Either $F_0 \cap E \not= \emptyset$ and $F_0 \cap B = \emptyset$, $\#F_0 \cap V = 1$ and $F_0 \cap B = \emptyset$ or $F_0 \cap B \not= \emptyset$. When $F_0 \cap B \not= \emptyset$ then do not consider. When $F_0 \cap E \not= \emptyset$ and $F_0 \cap B = \emptyset$ or $\#F_0 \cap V = 1$ and $F_0 \cap B = \emptyset$ then we have the above two cases. Therefore, if there is a face which is not in $D$ then we consider face in $D$ by the following way. Let $F_k$ be a face which is not in $D$ and $E(F_k) \cap B =\emptyset$. We get a sequence faces, say $F_{j1},$ $F_{j2},$ \dots, $F_{jr} = F_k$ where $F_{j1} \in D$ and $F_{js} \in F(N)\setminus D$ for $2 \le s \le r$. Since, $V(F_{js}) \cap V = \emptyset$ for $2 \le s \le r$, so, we follow the above argument for $F_{j2}$ and then, $F_{j3}$ and continue. At the end, we consider the face $F_k$. So, we follow the above argument for the faces whose have none of the edges in $B$ and none of the vertices in $V$. Similarly, if there is a face $F_x$ which contains edges of $B$ and has empty disjoint with $V$. In this case, we get a face sequence, say $F_{x1},$ $F_{x2},$ \dots, $F_{xt} = F_x$ where $F_{x1} \in D$ and $F_{xs} \in F(N)\setminus D$ for $2 \le s \le t$. Since $V(F_{xs}) \cap V = \emptyset$ for $2 \le s \le t$, so, we follow the above argument for $F_{x2}$ and then, $F_{x3}$ and continue. At the end, we consider the face $F_{x(t-1)}$. Hence the face $F_x$ will be in $D$ as $F_{x(t-1)}$ has a common edge with $F_x$. So, we follow this argument for the faces whose have some edges in $B$ and none of the vertices in $V$. This gives a tree $T$ which does not contain any vertices of $B$. Thus, we get a tree $T$ which satisfies the $2^{nd}$ and $3^{rd}$ property of proper tree. The set $D$ contains all the faces of $F(N)$. That is, $T$ touches all the faces of $M$. So, by lemma \ref{lem8}, $\#V(T) = n-2$. Therefore, $T$ is a proper tree. By theorem \ref{thm5}, $M$ contains a proper graph of type-II. So, by theorem \ref{thm2}, the map $K$ contains a contractible Hamiltonian cycle. Therefore, $K$ is Hamiltonian. This completes the proof of the theorem \ref{thm4}.
\hfill$\Box$

\section{Proper tree and proper graph of type-II \,:}\label{ptpg}

Let $M$ be a polyhedral map and $K$ denote its dual. We claim the following theorem.

\begin{theo}\label{thm5}  The edge graph $EG(K)$ contains a proper tree if and only if $EG(K)$ contains a proper graph of type-II.
\end{theo}

\noindent{\sc Proof of Theorem}\ref{thm5} : The map $K$ is dual of $M$. Let $T$ be a proper tree in $EG(K)$. Consider a graph $G(V, E)$ where $E = \{e \in EG(K)$ $|$ $V(e) \cap V(T) \not=\emptyset ~ \& ~ e \not\in E(T)\}$ and $V = \cup_{e \in E}V(e)$. We claim that $G$ is a proper graph of type-II. We first show that the dual of $G$ is a contractible Hamiltonian cycle. Then, by theorem \ref{thm2}, $G$ is a proper graph of type-II. Therefore, we only show that the dual of $G$ is a contractible Hamiltonian cycle. Let $T$ be dual of $2$-disk $D$ which is a tree. The boundary $\partial D$ is a Hamiltonian cycle and denote it by $C$. Let $e$ be an edge of $C$. We show that dual of $e$, say $e_d = uv$ belongs to $E$. The edge $e$ belongs to a face, say $F ( \in D)$. So, the dual vertex corresponding to $F$ belongs to $e_d$. That is, $V(e_d) \cap V(T) \not= \emptyset$. Let $F_u$ and $F_v$ be dual faces corresponding to $u$ and $v$ respectively. Since $e \in F$, so, either $F_u = F$ or $F_v = F$. Let $F_u = F$. By assumption $e = F_u\cap F_v$. If $F_v \in D$. Then, $e \not\in \partial |D|$. Thus, $F_v \not\in D$. That is, $e_d \not\in E(T)$. Again, $V(e_d) \cap V(T) \not= \emptyset$. Hence, $e_d \in G$. Therefore, dual edges corresponding to edges of $C$ belongs to $G$. Suppose, there is an edge $e1$ which belongs to $G$. Then, $V(e1) \cap V(T) \not=\emptyset $ and $e1 \not\in E(T)$. Suppose, it is not dual of any edge of $C$. By assumption, $V(e1) \cap V(T) \not=\emptyset$. So, the dual say, $e1_{d}$ belongs to interior of $D$. This implies that there are two faces in $D$ whose common edge is $e1$. So, by duality, the dual edge $e1_{d}$ belongs to the tree $T$. This is a contradiction as $e1 \not\in E(T)$. Therefore, the graph $G$ is the dual of $C$. Thus, the graph $G$ and $T$ define the same cycle. So, by theorem \ref{thm2}, $G$ is a proper graph of type-II.

Let $G_1(V_1, E_1)$ be a proper graph of type-II and $C_1$ be dual of $G_1$. By theorem \ref{thm2}, the cycle $C_1$ is contractible Hamiltonian cycle. It bounds a $2$-disk, say $D_1$. Consider dual of $D_1$ and denote it by $T_1$. We show that it is a proper tree. We follow the converse part of the {\em Theorem 2}\cite{maity upadhyay0}. It says that every contractible Hamiltonian cycle gives a proper tree in $K$. So, we consider the $2$-disk $D_1$ and follow the argument of {\em Theorem 2}\cite{maity upadhyay0}. This gives, the graph $T_1$ is a proper tree. This completes the proof of the theorem \ref{thm5}.
\hfill$\Box$

\begin{rem}\label{rem0} Let $T$ be a proper tree in $K$. Consider a graph $G :=\{e \in EG(K)$ $|$ $V(e) \cap V(T) \not= \emptyset \}$. We denote $G_1 := G \setminus T$. By the proof of theorem \ref{thm5}, the graph $G_1$ is a proper graph of type-II. Therefore, the graph $G$ can be decompose into proper tree and proper graph of type-II. That is, $G  = G_1 \cup T$.
\end{rem}

\section{The steps for searching Hamiltonian cycle in equivelar maps \,: }\label{algorithm}
The following steps may be implemented as a computer program to search contractible, non-separating and noncontractible separating Hamiltonian cycles. We use definition of face sequence and some notations of lemma \ref{lem4} in the following algorithms.

\medskip

\begin{algo}\label{algo1}  Let $EG(K_1)$ be the edge graph of a $\{p, q\}$ equivelar map $K_1$ and $\#V(EG(K_1))$ $= n$. Let the set $K$ contains faces of the polyhedral map $K_1$. We follow the following steps.
\begin{enumerate}
\item Construct dual of $K$. It is a set say, $M := \{F_1, \dots, F_n\}$ of faces where $|M|$ is dual of $K_1$.
\item Put all edges of $M$ in $E_M$.
\item The number of vertices of $K$ is $n$. That is, $\# V(K) = n$. Let $S_1 = \emptyset$. Then, consider $n$ edges from $E_M$ and put in $S_1$. We follow the following steps.
  \begin{enumerate}
   \item[A.] If $\# F \cap S_1 = 2$ for all $F \in M$ then we continue to next step. Otherwise, we consider another possible $n$ edges from $E_M$ and continue.
   \item[B.] Let $S_2 = \emptyset$. Put, dual edges of $S_1$ in $S_2$. Now, there are two possibilities.
       \begin{enumerate}
       \item[a.] Check existence of a face sequence $FS(F_{t_1}, F_{t_n})$ in $M$ such that $CEFS(F_{t_1},$ $F_{t_n})$ $= S_1$. Suppose, such face sequence exists. By the argument of theorem \ref{thm1}, the set $S_2$ is a cycle.

       \item[b.] Suppose, there is no such sequence $FS(F_{t_1}, F_{t_n})$. We have $F \cap S_1$ $\forall$ $F \in M$. So, we get disjoint collections of faces, namely, $H_1, \dots, H_r$ in $M$ where each $H_i$ has above property $(a)$. That is, let $H_i = \{F_{i_1}$, \dots, $F_{i_m}\}$. Then, $F_{i_1} \cap F_{i_m} \in S_1$ and  $F_{i_j} \cap F_{i_{j+1}} \in S_1$ for $1 \le j \le m-1$. Let $e_i = F_{i_1} \cap F_{i_m}$ and $e_j = F_{i_j} \cap F_{i_{j+1}}$ for $1 \le j \le m-1$. Consider dual edges corresponding to $e_i$ and put in $T_i$. The set $\{e_1, \dots, e_m\}$ satisfies the $2^{nd}$ and $3^{rd}$ property of Definition \ref{defn1}. So, $T_i$ is a cycle. Thus, $S_2$ contains $r$ disjoint cycles.
       \end{enumerate}
   \item[C.] If $r=1$ then, we consider $S_2$. This is Hamiltonian cycle. If $r > 1$ then we go to step $(3)$. Choose another possible $n$ edges from $E_M$ and continue.
  \end{enumerate}
 \item When $S_2$ is a Hamiltonian cycle then we go to the following steps to classify the cycle.
 \begin{enumerate}
 \item[(i)] Suppose, $\frac{n - 2}{p -2}$ is not an integer. We use the result of \cite{maity upadhyay0}. Thus, $S_2$ is non-contractible Hamiltonian cycle. To classify further, consider the set $G(V_M, E_M \setminus S_1)$. We have the following two cases.
     \begin{enumerate}
     \item[a.] If $G(V_M, E_M \setminus S_1)$ is connected, then, the cycle $S_2$ is separating Hamiltonian cycle (by theorem \ref{thm1}).

     \item[b.] If $G(V_M, E_M \setminus S_1)$ has two components and none of them is proper tree, then, the cycle $S_2$ is noncontractible separating Hamiltonian cycle (by theorem \ref{thm3}).
      \end{enumerate}
 \item[(ii)] Suppose $\frac{n - 2}{p -2}$ is an integer. We check graph $G(V_M, E_M \setminus S_1)$. We go to the following cases.
    \begin{enumerate}
     \item[a.] If $G(V_M, E_M \setminus S_1)$ is connected, then, the cycle $S_2$ is separating Hamiltonian cycle (by theorem \ref{thm1}).
     \item[b.]  If $G(V_M, E_M \setminus S_1)$ has two components and one of them is proper tree, then, the cycle $S_2$ is contractible Hamiltonian cycle (by theorem \ref{thm2}).
     \item[c.] If $G(V_M, E_M \setminus S_1)$ has two components and none of them is proper tree, then, the cycle $S_2$ is noncontractible separating Hamiltonian cycle (by theorem \ref{thm3}).
      \end{enumerate}
 \end{enumerate}
 \item We stop unless we get contractible, non-separating and noncontractible separating Hamiltonian cycle (if exists) after finite steps.
\end{enumerate}

\end{algo}

\section{The steps for searching a non-contractible Hamiltonian cycle in equivelar maps \,: }\label{algorithm1}

The following steps (using backtracking) may be implemented as a computer program to locate non-contractible Hamiltonian cycles\,:
\medskip

\begin{algo}\label{algo2} Let $EG(K)$ be the edge graph of a $\{p, q\}$ equivelar map $K$ and $\#V(EG(K)) = n$. Let $M$ be the set of all  $p$-gonal faces and $i$ denote the number of steps. We construct two set $D$ and $V$ as follows : choose an element $P_{0}\in M$. Define $D:=\{P_{0}\}$, $V := \{V(P_{0})\}$ and $i = 1$. We have either $\#V = \#V(EG(K))$ or $\#V < \#V(EG(K))$. We go to the following steps.
\begin{enumerate}

\item If $\#V = \#V(EG(K))$. We go to the step $(3)$.
\item If $\#V < \#V(EG(K))$, then, we observe the set $D$ at $i^{th}$ and $(i+1)^{th}$ steps. Let $v\in V(K)\setminus V$. By the {\em Algorithm 1}\cite{maity upadhyay0}, there exists a face $P$ such that $v \in V(P)$, $V(P)\cap V = \{v_{1},v_{2}\}$ and $E(P)\cap E(P_1)=\{\{v_{1},v_{2}\}\}$ for some $P_1\in D$. Then put $D = D\cup \{P\}$, $V = V\cup V(P)$ and $i = i + 1$. We check the sets $D$ and $D \cup \{P\}$. If $D \cup \{P\}$ is a $2$-disk then we put $D = D \cup \{P\}$ and $V = V\cup V(P)$ and $i = i+ 1$. Go to the next step and continue. Do this until we get either $V = V(EG(K))$ or $D$ is a $2$-disk and $D \cup \{P\}$ is not a $2$-disk. Then, we go to next step $(3)$.

\item If $V = V(EG(K))$, then, we add a face $P_2$ of $M \setminus D$ in $D$. The geometric carrier $\partial |D \cup P_2|$ contains more than one cycle as $V(P_2) \subset V$. Denote it by $C_1$. Otherwise, at the end of processes in step $(2)$, we get a face $F$ whose at least one vertex is already in $V$. In this case, we consider the geometric carrier $\partial | D \cup F|$ which contains at least one cycle and denote it by $C_1$. Here, at both the cases, the cycle $C_1$ does not bound any $2$-disk. So, the cycle is non-contractible. Put set $V_1 := V(C_1)$. We go to next step.

\begin{enumerate}
\item If $\#V_1 = \#V(EG(K))$. We stop here. We get a non-contractible Hamiltonian cycle $C_1$
\item If $\#V_1 < \#V(EG(K))$, then, there is a vertex in $V(K) \setminus V$. We follow the step $(2)$. We get a new face $P_3$ at $j^{th}$ step. The cycle $C_1$ and $\partial P_3$ have common path. So, we concatenate $C_1$ and $\partial P_3$ if the concatenate cycle has length bigger than the length of $C_1$. We get new cycle. Denote it by $C_1$ and put, $V = V(C_1)$ and $j = j+ 1$. Go to next step and continue. Do this until we get $V_1 = V(EG(K))$. This gives a non-contractible Hamiltonian cycle $C_1$.

\end{enumerate}

\end{enumerate}

  We stop unless we get a non-contractible Hamiltonian cycle (if exists) after finite steps.

\end{algo}

\section{Acknowledgement}

Work of second author is partially supported by SERB, DST grant No. SR/S4/MS:717/10.

\end{document}